\documentclass[11pt]{amsart}
\usepackage{geometry}                
\usepackage{graphicx}
\usepackage{amssymb}
\usepackage{color}
\definecolor{Yellow}{rgb}{1,0,0}
\definecolor{Green}{rgb}{0,1,0}
\usepackage{amsmath}
\newcommand{\cdo}{}
\newcommand{\rjm}{} 
\newcommand{\rjo}{}
\newcommand\dsum{\sum}
\newcommand\ma{\mu}

\newcommand{\diag}{\text{diag}}

\newcommand{\p}{\partial}

\newtheorem{theorem}{Theorem}

\newtheorem{claim}[theorem]{Claim}

\newtheorem{corollary}[theorem]{Corollary}

\newtheorem{definition}[theorem]{Definition}

\newtheorem{lemma}[theorem]{Lemma}

\newtheorem{remark}{Remark}

\title[Unique equilibria and substitution effects in the marriage market]
{Unique equilibria and substitution effects in a stochastic model of the marriage market}

\begin{document}

\thanks{
The authors are grateful to Aloysius Siow for attracting
their attention to this question, and for many fruitful discussions.
This project developed in part from the Master's research of CD; however, the present
manuscript is based a new approach which greatly extends (and largely subsumes)
the results of his thesis  \cite{Decker10},
\rjm{and of the 2009 preprint by three of us entitled
{\em When do systematic gains uniquely
determine the number of marriages between different types in the Choo-Siow
matching model? Sufficient conditions for a unique equilibrium.}
}
CD, RJM and BKS are pleased to acknowledge the support of Natural Sciences and
Engineering Research Council of Canada (NSERC) grant 217006-03 RGPIN; CD also
benefitted from an Undergraduate Student Research Award (URSA) held in the
summer of 2009. EHL acknowledges a United States National Science Foundation
grant PHY-0965859. \copyright 2011 by the authors.}
\author{Colin Decker}
\address{Department of Mathematics, University of Toronto, Toronto
Ontario Canada M5S 2E4 \texttt{colin.decker@utoronto.ca}
Current address: SunLife Financial, 150 King Street West Toronto Ontario
Canada M5H 3T9 ATTN:TK11 Corporate Risk Management
}
\author{Elliott H. Lieb}
\address{Departments of Mathematics and Physics, Jadwin Hall, Princeton
University, P.O. Box 708, Princeton, NJ 08542, USA
\texttt{lieb@math.princeton.edu}}
\author{Robert J. McCann}
\address{Corresponding
author: Department of Mathematics, University of Toronto, Toronto Ontario
Canada M5S 2E4 \texttt{mccann@math.toronto.edu}}
\author{Benjamin K. Stephens}
\address{Department of Mathematics, University of Toronto, Toronto
Ontario Canada M5S 2E4 \texttt{stephens@math.toronto.edu}. Current address:
Department of Mathematics, University of Washington, Seattle Washington USA
98195-4350 \texttt{benstph@math.washington.edu}}

\date{\today}

\begin{abstract}
Choo-Siow (2006) proposed a \rjo{model for the marriage market which allows for random}
identically distributed noise in the preferences of each of the participants. The
randomness is
McFadden-type, which permits an explicit resolution of the equilibrium preference probabilities.

\rjo{The purpose of this note is to prove} 
uniqueness of the resulting equilibrium marriage distribution,
and find a representation of it in closed form.  This allows us to derive
smooth dependence of this distribution on exogenous
preference and population parameters, and establish sign, symmetry, and 
size
of the various substitution effects, \rjo{facilitating} comparative statics.  For example,
we show that an increase in
the population of men of any given type in this model
leads to an increase in single men of each type,  and a decrease in single women of each type.
We show that an increase in the number of men of a given type increases the equilibrium transfer
paid by such men to their spouses, and also increases the percentage of men of that type who
choose to remain unmarried.  While the above trends may not seem surprising,
the verification of such properties helps to substantiate the validity of the model.
Moreover,  we make unexpected predictions
which could be tested:
namely,  the percentage change
of type $i$ unmarrieds with respect to fluctuations in the total number of
type $j$ men or women turns out to form a symmetric positive-definite matrix
$r_{ij}=r_{ji}$ in this model, and thus to satisfy bounds
such as $|r_{ij}| \le (r_{ii}r_{jj})^{1/2}$.

\rjo{Along the way, we give a new proof for the existence of
an equilibrium, based on a} \rjo{strictly convex}
variational principle and a simple estimate, rather than a fixed point theorem.
Fixed point approaches to the existence part of our result have been explored by others
\cite{CSS} \cite{Dag} \cite{Fox},  but are much more complicated and yield neither uniqueness, nor
comparative statics, nor an explicit representation of the solution.

JEL Classification: J12, C62, C78, C81, D03, Z13

KEYWORDS: Choo-Siow, marriage market,
matching, random, unique equilibrium, comparative statics, convex analysis.


\end{abstract}

\maketitle

\section{Introduction}
\label{S:introduction}

\rjo{The classic transferable utility framework which Becker used to model the marriage market
was augmented by Choo and Siow  \cite{CS} to allow for the possibility that
agents' preferences might be only partly determined by observable characteristics,
and might therefore include a stochastic component
depending on unobservable characteristics.}
The randomness was chosen to be
McFadden type \cite{McFadden}, and spreads the preferences of agents on one side of the
marriage market over the entire type-distribution of agents on the other,
\cdo{thus yielding
non-assortative matching, a feature of the marriage market that has long been observed empirically \cite{CS}. The model is non-parametric
and highly tractable; given an observation of marriages between different types of agents, there is a simple closed form, \rjo{`}point-identified' expression for the average
utility (or \rjo{`}total gains') generated by each type of marriage. Of use to econometricians, demographers, and economic theorists is a solution to the backwards problem:
given the total gains to each type of marriage and population data, what can be said about
\rjo{possible distributions of marriages corresponding to those parameters?
Basic questions include existence and uniqueness of such a
marital distribution; computing comparative statics is even more interesting}.
How will the marital distribution respond to
changes in composition of the population? How will it change due to policy shocks
\rjo{whose effects can be summarized by revising}
the total gains parameters?
The purpose of the present note is to address these questions.}

A priori, it is not evident whether each choice of
gains parameters leads to the existence of a
distribution of marriages which clears the market
 \cite{CSS} \cite{Fox} \cite{Dag}.
\rjo{Although it is not the main point of our paper, we begin by reconfirming that it
does.
We then go on to} show this equilibrium is unique,  and give an explicit formula
for the resulting marriage distribution, essentially solving the model completely.
This allows us to derive smooth dependence of the resulting equilibrium
state on the specified population and preference parameters,  and establish sign,
symmetry and bounds on the response of the predicted distribution of marriages
to changes in each parameter.
Although existence of an equilibrium was also discussed by Choo, Seitz and Siow \cite{CSS}
and Fox \cite{Fox} (and by Dagsvik \cite{Dag} for a related model),
our proof relies on 
the reformulation of the problem as a variational minimization,
hence is simpler than the fixed point argument they suggest.
More \rjo{importantly}, 
it yields the solution in closed form.
Our uniqueness result
is the first concerning this model
\cite{CSW},
and is based on convexity (in appropriate variables) of the new variational principle formulated
in Section \S\ref{S:Thm 1}.
The variational technique we use is powerful, because it characterizes the
endogenous variables of interest as the critical points of a function.
Indeed, the competitive equilibrium in this model turns out to be realized
as the \rjo{minimum} of a strictly convex function (given in Theorem \ref{T:main} below),
\rjo{facilitating its subsequent analysis.}

\subsection{{\rjo{Organization}}}

\rjo{The remainder of this introduction provides further motivation and background
for the Choo-Siow model.
Section \S\ref{S:CS model} details the model and could be skipped by readers
familiar with the Choo-Siow model,  except that it ends with a
summary placing our results in the context of related literature.
Further comments concerning the derivation of the model
may be found in Appendix \ref{A:preference probability}.}
 Section \S\ref{S:results}
states our results formally; it is followed by a section containing
derived statics and one further conjecture.
Section \S\ref{S:Thm 1} proves the existence and uniqueness
of equilibria, while \S\ref{S:Thm 2} is devoted to the comparative statics
established in Theorem \ref{T:comparative statics}.
Section \S\ref{S:comparative statics} establishes the further comparative static
assertions of Corollary \ref{C:transfer and percent} and \S \ref{S:derived conjectures}.

\subsection{{\rjo{Further remarks}}}
The random component of agent preferences is a salient feature of the Choo-Siow model.
Due to this randomness,
the equilibrium marriage distribution predicted by the model
will not be positive assortative, even when the observed
attributes of the agents are one-dimensional. This is consistent with empirical data.
Even in experiments where the agents are parameterized by ordered types,
such as age,  observed marital data will almost never be genuinely assortative.  For example,
in any given population it is unlikely to be true that the age of the youngest woman married
to a 34 year old man always exceeds that of the oldest woman married to a 33 year old.
Similarly, one always finds matches in large populations that pair high with
low qualities as measured by any given ordered observable characteristic (e.g. income, years of education, etc.).
Thus a strictly assortative framework fails to explain the presence of, for example,
the existence of PhD graduates married to high-school drop-outs.

The classic transferable utility model of the marriage market, introduced by
Gary Becker  \cite{Becker}, in principle predicts how agents will marry given
exogenous preference parameters. However, it has seldom been estimated. There are two main
obstacles \rjo{to} estimating a model of the marriage market. First, equilibrium transfers
in modern marriages (except in the case of dowries) are not observed. Hence any
behavioural model that requires their presence in data is not identifiable. Second,
real-world agents are described by discrete, multi-dimensional, possibly unordered,
types. But the classic Becker model predicts positive assortative matching under
the assumption that agent type is one dimensional, continuous, ordered, and that
preferences are super-modular. This positive assortative matching is limiting, but
does ensure that his model predicts a unique marital distribution.

The Choo-Siow model eliminates the structural assumptions of the classic model. First, it is not necessary to observe transfers in order to determine the equilibrium marriage distribution generated by the model. In fact,
we provide an explicit formula for the equilibrium marriage distribution in terms of the
derivative of the Legendre transform of a known function. Second, the model places no a priori
structure on the nature or number of types that agents (men and women) can have. This allows
consideration of a wide range of attributes, like race, religion, level of income, and educational
achievements.

In this more realistic framework, with its lack of structure for
the agents' deterministic preferences and types, the issue of whether there exists an
equilibrium marital distribution, and if so whether it is unique, becomes a question
of fundamental theoretical and econometric significance.
The theoretical importance arises from the fact that 
uniqueness of equilibria in two-sided matching problems is usually not better than
a generic property,
except perhaps in certain convex programming settings like \cite{Ekeland10} \cite{FigalliKimMcCann-econ},
which include continuous Monge-Kantorovich matching \cite{GOZ} \cite{ChiapporiMcCannNesheim10}.
Further, the randomness considered in the model below is the commonly used
extreme value logit type, thus any result that describes properties of the
equilibrating matches has potentially wider applicability. The econometric importance arises
from the fact that models of the marriage market are useful to econometricians only insofar as
they make unique predictions of a marital distribution, given exogenous preferences.
From a practical point of view having closed form solutions which permit
comparative statics may be even more crucial.

\section{The Choo-Siow Marriage Matching Model}
\label{S:CS model}

Our presentation emphasizes the stochastic heterogeneity that differentiates the Choo-Siow model
from classical models. The competitive framework, which uses transfers of utility from spouses
to equilibrate the market, is explored in detail in Choo-Siow \cite{CS} but treated here
only in \S\ref{S:comparative statics}. 
It should be noted at the outset  that the methods developed
here also apply to other non-transferable utility models present in the literature. For example,
Dagsvik \cite{Dag} develops a model of the marriage market which uses an assignment algorithm
(deferred acceptance) rather than utility transfers to sort matches, but his equilibrium
conditions are functionally similar to ours.

\subsection{Setting}

What is exogenous in this model are the observed types of men
and of women,  the numbers of men and women of each type in the population,  and the total gains
$\pi_{ij}$ of marriage between a man of observed type $i$ and a woman of
observed type $j$, relative to both partners remaining single.  The quantity
$\pi_{ij}$ will not reappear until \eqref{sum}.
On the other hand, individual agents have a utility function that depends on both an endogenous
deterministic component that captures systematic utility, and an exogenous random one that
models heterogeneity within the population of each given type.
Thus the utility accrued by a man of type $i$ and specific identity $g$ who marries a woman
of type $j$ is assumed to be:
\begin{equation}
V_{ijg}^{m}= \eta^{m}_{ij}+\sigma \epsilon_{ijg};
\end{equation}
the case $j=0$ represents the utility of remaining single.
The deterministic component is $\eta^m_{ij}$; its endogeneity can be interpreted
to reflect the possibility of interspousal transfer, as in Choo-Siow  \cite{CS} and
\S\ref{S:comparative statics} below.
It is set in equilibrium, and depends explicitly on the type of the man and the type of the
woman, and implicitly on market conditions, i.e. on the relative abundance or scarcity of men
and women of each different type.
The random term $\epsilon_{ijg}$ depends additionally on the specific identity of the man,
but not on the specific identity of the woman.  Hence a specific thirty-five year old man may
have stronger than typical (with respect to his age group) attraction for fifty-year old women.
But this attraction does not depend on whether, for example, the older woman has an
especially strong attraction to younger men
(assuming this latter characteristic is unobservable in the data and hence not reflected in $j$).

The random term is assumed to have the Gumbel extreme value distribution
described in Appendix \ref{A:preference probability}.
This distribution was introduced to the economics literature by McFadden \cite{McFadden}.
\rjo{It represents a severe simplifying assumption in our model,  but has recently relaxed
in work by Galichon and Salani\'e~\cite{GalSal2} which allows heteroskedasticity, for example.
Our homoskedastic model admits
only a single scaling parameter $\sigma$ which} measures the degree of randomness; its reciprocal
can be interpreted as the signal to noise ratio. It is equal to unity in the original
Choo-Siow model. For illustrative purposes, we will have occasion to allow $\sigma$ to vary
and in doing so embed the Choo-Siow model in a one parameter family of models that differ by
the degree of randomness present in them.

In contrast with deterministic matching models, agents of a particular type do not have a uniform
preferred match. This is because their preferences depend on the random variable $\epsilon_{ijg}$.
Using the Gumbel structure, the probability that a man of type $i$ prefers a woman of type
$j$ among all other possible marital choices $k \in \{0, ..., J\}$ is given
by
\begin{equation}\label{probability}
\mbox{Pr(Man of type $(i,g)$ prefers a woman of type $j$)}=
\frac{\exp(\frac{\eta^{m}_{ij}}{\sigma})}{\sum_{k=0}^{J}\exp_(\frac{\eta^{m}_{ik}}{\sigma})};
\end{equation}
(see Appendix \ref{A:preference probability} for a derivation).
This probability distribution is endogenous, because it depends on the various $\eta^m_{ik}$.
Note that it does not depend on the specific identity $g$ of the man of type $i$,
since the noise is identically distributed for each different $g$.
Yet it is possible already to see how the equilibrium marriage output will differ markedly
from a deterministic one. Whereas in the deterministic case all members of a given type
typically have the same preferred match, here the preferred matches of type $i$ men are smeared across all
female types according to the distribution defined by (\ref{probability}). The mean and spread
of the smearing are determined by the endogenous values $\eta^{m}_{ij}$ and by $\sigma$,
respectively.

Consider $\sigma \in [0,\infty]$. The case $\sigma = 1$ corresponds to the Choo-Siow model where
some smearing is present. The case $\sigma = 0$ corresponds to a deterministic matching model,
for which there is no smearing. Indeed, as $\sigma \rightarrow 0$, the largest exponentials
dominate all others, and the probability that a man of type $i$ prefers a woman of type $j$
converges to $0$ or $\frac{1}{\# { \arg\max \{\eta^{m}_{ij}\mid 0 \le j \le J\} }}$ depending on
whether or not $\eta^{m}_{ij}$ weakly dominates all other preference parameters $\eta^{m}_{ik}$.
Conversely, as $\sigma \rightarrow \infty$, the stochastic term dominates the utility function, and the resulting probability distribution converges to the uniform distribution. In this case, there is maximal smearing,
as preferences are completely random, constrained only by availability of prospective partners
to marry.

Female preferences are also smeared, and the equilibrium marriage distribution is determined when $\eta^{m}_{ij}$ and $\eta^{f}_{ij}$ are such that the number of desired marriages of each type is the same on both sides of the market.

We now elaborate on the Choo-Siow model.  We henceforth fix $\sigma =1$;
since preferences are relative, this normalization
can always be attained by rescaling all of the preferences in the model.


\subsection{The Choo-Siow model}
\label{SS:CS model}

Suppose we wish to predict the number of marriages between men and women of different types.
The number of men of type $i$ is denoted  $m_{i}$. The number of marriages of
type $i$ men to type $j$ women is denoted $\mu_{ij}$. \rjo{The number of type $i$ men
(respectively type $j$ women) who choose to remain single is denoted by
 $\mu_{i0}$ (and $\mu_{0j}$
respectively).}
If each man marries his preferred woman, the equality
\begin{equation}\label{LLN}
\mbox{Pr(Man of type $(i,g)$ prefers a woman of type $j$)}=\frac{\mu_{ij}}{m_{i}}
\end{equation}
will be valid, or at least 
as the population size becomes large, the right hand side of the equality converges
to the left hand side by the law of large numbers, or the maximum likelihood theorem.

Using equations \eqref{probability}--\eqref{LLN} to compute the ratio of the probability
that a man of type $i$ prefers a woman of type $j$ to the probability that he prefers
to remain single, we arrive at the following formula:
\begin{equation}{\label{utility}}
\mu^{m}_{ij}=\frac{e^{\eta^{m}_{ij}}}{e^{\eta^{m}_{i0}}}\mu^{m}_{i0}.
\end{equation}
These $I \times J$ equations are in fact quasi-demand equations, because they indicate the number of type $\mu_{ij}$ marriages that men of type $i$ would like to participate in. Viewing the female market cohort as the supply side, there are analogous supply equations.
Letting the utility acquired by a woman of type $j$ and specific identity $h$ who marries
a man of type $i$ be
\begin{equation}
V^{f}_{ijh}= \eta^{f}_{ij}+\epsilon_{ijh},
\end{equation}
the above analysis produces $I \times J$ supply equations of the form:
\begin{equation}\label{pre-equilibrium}
\mu^{f}_{ij}=\frac{e^{\eta^{f}_{ij}}}{e^{\eta^{f}_{0j}}}\mu^{f}_{0j}.
\end{equation}
The equilibrium output in the Choo-Siow model is a specification of $\mu_{ij}$ for all
$0\leq i \leq I$, and all $0 \leq j \leq J$. This output is obtained by requiring that
supply balance demand: $\mu^{f}_{ij} = \mu^{m}_{ij}$. Under this market-clearing hypothesis,
we have the following equation. The endogenous parts of the
$\eta^{m}_{ij}$, $\eta^{f}_{ij}$, $\eta^{m}_{i0}$, $\eta^{f}_{0j}$ are eliminated upon adding
them to arrive at the definition of an
exogenous\footnote{See  \S\ref{S:comparative statics} below or Choo and Siow \cite{CS} for an explanation in terms of spousal transfers.}
aggregated gains variable $\pi_{ij}$ associated to each observable type of marriage:

\begin{equation}\label{sum}
\pi_{ij} := \frac{\eta^{m}_{ij}+\eta^{f}_{ij}-\eta^{m}_{i0}-\eta^{f}_{0j}}{2}.
\end{equation}
Using the market-clearing hypothesis, we may re-write the equilibrium condition
$\mu^m_{ij} = \mu^f_{ij} =:\mu_{ij}$ in terms of the exogenous variable $\pi_{ij}$ as follows:

\begin{equation}
\frac{\mu_{ij}}{\sqrt{\mu_{i0}\mu_{0j}}}=e^{\pi_{ij}}.
\end{equation}
Finally, letting $\Pi_{ij} =e^{\pi_{ij}}$, the equilibrium output is given by
\begin{equation}\label{equilibrium}
\frac{\mu_{ij}}{\sqrt{\mu_{i0}\mu_{0j}}}=\Pi_{ij}.
\end{equation}
The equilibrium conditions expressed in equation (\ref{equilibrium}) are implicit.
They give necessary conditions for real numbers $\mu_{ij}$ to be an 
output of the Choo-Siow model. However, they are not sufficient; a secondary
set of necessary conditions, population constraints, must also be satisfied.
Let there be $I$ types of men, and $J$ types of women. The number of men of type $i$
is denoted $m_{i}$, and the number of women of type $j$ is denoted $f_{j}$.
The vector whose $i^{th}$ component is $m_{i}$, and whose $(I+j)^{th}$
component is $f_{j}$, is denoted by $\nu $. 
Called the \textit{population vector},
it has $(I+J)$ components and may also be denoted by $[\,m\mid f\,]$.
A specification $(\ma_{ij})_{0 \le i \le I \atop 0 \le j \le J}$
of the number of marriages and
singles of each type
is called a \textit{marital distribution.}
The following population
constraints must be satisfied by all marital distributions, and are a
consequence of the definitions:\qquad \qquad

\begin{eqnarray}\label{males}
\ma_{i0}+\dsum_{j=1}^{J}\ma_{ij} &=&m_{i}, \\
\label{females}
\ma_{0j}+\dsum_{i=1}^{I}\ma_{ij} &=&f_{j}, \\
\label{non-negative}
\ma_{ij}\ge 0.
\end{eqnarray}
\cdo{Several questions naturally arise.
We call these questions the \textit{Choo-Siow inverse problem}:}

\smallskip\noindent
\textbf{Problem (Choo-Siow inverse problem)} \textit{Given a gains
matrix $\Pi= (\Pi_{ij})$ and a population vector $\nu= [\,m \mid f\,]$, does there
exist
a unique marital arrangement generating $\Pi$? In other words, assuming
the entries $\Pi_{ij}$ to be non-negative and $m_{i}$ and $f_{j}$ to be
strictly positive, does exactly one matrix $(\mu_{ij})$ with non-negative
entries exist\footnote{\rjo{Some readers have pointed out that for each realization of the
randomness,
balancing supply with demand amounts to solving a linear program.
In this setting,  existence of an equilibrium (and its generic uniqueness) are
well-known \cite{GOZ} \cite{ChiapporiMcCannNesheim10}.  However,  instead of balancing supply
with demand for each realization,  the Choo-Siow model is based on balancing
{\em expected} supply \eqref{utility} with {\em expected} demand \eqref{pre-equilibrium}.
Since the expected supply
(or demand) need not correspond to the actual supply (or demand) for any
realization of the randomness, existence of an equilibrium in the sense prescribed by
Choo and Siow is not obvious.  Unlike a Nash equilibrium,  in which no individuals
have both the incentive and the opportunity to change their marital status,  in a
Choo-Siow equilibrium this will only be true in some average sense.  Never the less,
the notion has proved useful empirically,  and might even be construed to
reflect the metastability (as opposed to stability) displayed by actual marriage markets.
As simple examples show, such equilibria generally correspond to interior
points in the feasible polytope defined by the population constraints,  hence cannot
be selected by the minimization of any {\em linear} function defined on this polytope.
Nevertheless, we shall show they do exist and are selected by the minimization of a
specific strictly convex function introduced below.}}
that satisfies \eqref{equilibrium}--\eqref{non-negative}?
\rjo{Furthermore, can the qualitative dependence of $(\mu_{ij})$ on the exogenous parameters
$\Pi=(\Pi_{ij})$ and $\nu = [\nu_i]$ be described?}}

\smallskip
This problem is important for several reasons.
First, the implicit conditions present in equation (\ref{equilibrium}) are the equilibrium outcome of a
competitive market. There are not so many realistic environments with
finitely many agent types and many commodities which are known to generate
unique competitive equilibria --- except possibly generically.
While there are generic uniqueness results for matching problems
that can be reduced to convex programing problems such as Monge-Kantorovich matching,
e.g.\ \cite{GOZ} \cite{ChiapporiMcCannNesheim10} \cite{Ekeland10} \cite{FigalliKimMcCann-econ},
the stochastic heterogeneity prevents the equilibrium in our model from being formulated as such.
\rjo{Instead, stochasticity restores uniqueness without the need for a genericity assumption
in our model --- and indeed in the more general setting studied by
Galichon and Salani\'e \cite{GalSal2}.}

Second, an affirmative explicit solution to the \textit{Choo-Siow Inverse Problem} makes the
Choo-Siow model useful in econometric analysis. The matrix $\Pi$ is exogenous and unobserved
in data, but can be point-estimated from an observed marriage distribution. An economic or
social shock will affect the systematic utilities that agents of various types incur by
marrying agents of various others, and will therefore alter the value of $\Pi$. This effect
can be approximated to form an updated matrix of aggregated systematic parameters
$\Pi^{'}$. Existence and uniqueness guarantee that there will be exactly one marriage
distribution that results from the shock, making the model predictive. In the same vein,
demographers are often interested in predicting how marriage distributions will change
due to changing demographics, i.e. changes in the population vector $\nu$.
Our closed form solution makes it possible to compute the sign and in some cases the
magnitude of such changes explicitly.

Finally, if the model can be shown to admit a unique distribution, the estimated parameters
$\pi_{ij}$ are an alternative characterization of the observed marriage distribution. The
recharacterization is useful because the parameters of the Choo-Siow model have a behavioural
interpretation, and are not merely observed data.

\subsection{Summary of progress and related literature}
\label{S:progress}

\rjo{The} related {\em local} uniqueness question was resolved
by Choo and Siow in \cite{CS}.
However the issue of {\em global} uniqueness was left open, and posed as an open
problem in a subsequent working paper
by Siow \cite{Siow08}.  We resolve this question positively by introducing a variational principle and a change of variables which allows us to exploit convexity.
The question of existence of $(\mu_{ij})$ for all $\Pi = (e^{\pi_{ij}})$ was
addressed 
in a working paper
of Choo, Seitz and Siow \cite{CSS} by appealing to the Tarski fixed point theorem;
see also the related results of Fox \cite{Fox} and Dagsvik \cite{Dag}.
However the proofs there
are long and \rjm{involved,} 
whereas the variational proof in the present paper
is simple and direct and follows from continuity and compactness by way of an elementary estimate.
Moreover, it leads to an explicit representation of the solution.  This allows us to rigorously
confirm various desirable and intuitive features of Choo-Siow matching,  whose presence or absence
might in principle be used as a test to refute the validity of various alternative matching models.
Among other results, we show for example that an increase in the number of men of a given type
increases the equilibrium transfer paid by such men to their spouses,  while also increasing the
percentage of such men who choose to remain single.  See Theorem \ref{T:comparative statics}
below for related statements and more surprising conclusions.
%
%
\rjm{Independently of us,
Galichon and Salani\'e \cite{GalSal} \cite{GalSal2}
developed a variational approach which extends
the Choo-Siow model by allowing for randomness of much more general form.
From theoretical considerations,  they derive a
closed-form expression for a
strictly concave social welfare function governing competitive equilibria in their
models, and identities relating different sets of endogenous variables under consideration.
This gives an alternate approach to existence and uniqueness of equilibria in
our setting as a special case. 
Apart from its convex analytic nature, their approach is quite distinct
from ours, and 
leads to very different insights and interpretations.}


\section{Precise Statement of 
Results}
\label{S:results}

In the preceding remarks, the \textit{Choo-Siow inverse problem} was phrased in terms of finding
existence and uniqueness of equilibrium $\mu$ given exogenous data $\Pi$ and $\nu$.
As the name suggests, it is also useful to think of this problem as one of inverting a function.
From this point of view, even though $\Pi$ is exogenous, we may prefer instead to consider $\Pi$
as the image of a marriage distribution under a transformation that we seek to invert.
\begin{remark}[Incomplete participation]
From the market equilibrium point of view, the fact that the left hand-side of
\eqref{equilibrium} becomes infinite when $\mu_{i0}$ or $\mu_{0j}$ is equal to zero
is unproblematic. It means that for no finite value of the exogenous $\Pi$ is sufficient
to induce all the representatives of some type to marry.
However from the inverse
problem point of view, it is necessary to stipulate that $\mu_{i0}$ and $\mu_{0j}$ be
strictly positive.
\end{remark}



\subsection{Preliminaries}

Let us begin with a reformulation of the Choo-Siow inverse problem; Siow 
attributes this reformulation to Angelo Melino.
Let $\beta_i := \mu_{i0}$ and $\beta_{I+j}:=\mu_{0j}$ denote the number of unmarried
men and women of types $i=1,\ldots,I$ and $j=1,\ldots,J$ respectively.
Since the gains matrix \eqref{equilibrium} can be used to express each component
$\mu_{ij} = \beta_{i}\beta_{I+j}\Pi_{ij}$ of the marital distribution
in terms of these new variables, the population constraints
\eqref{males}--\eqref{females} can be reduced to a system
%

\begin{eqnarray}
\label{ast}%
\beta_{i}^{2}+\sum_{j=1}^{J}\beta_{i}\beta_{I+j}\Pi_{ij}-\nu_{i}=0, && 1\leq
i\leq I,\cr
\beta_{I+j}^{2}+\sum_{i=1}^{I}\beta_{i}\beta_{I+j}\Pi_{ij}-\nu_{I+j}=0, && 1\leq
j\leq J.
\end{eqnarray}
of $(I +J)$
quadratic polynomials in the $(I+J)$ variables $\{\beta_{k}\}_{k=1}$ counting
the number of umarried men and women of each type.

A solution to this system of equations is a vector of amplitudes $\beta$ that
has $(I + J)$ components. Abstractly, its components might be real, complex,
or both. The Choo-Siow Inverse Problem is equivalent to showing that the
polynomial system (\ref{ast}) has a unique solution with real positive
amplitudes
for all gains matrices $\Pi$ 
and population vectors $\nu=[\,m \mid f\,]$ with positive components.
The full marital distribution satisfying \eqref{equilibrium}--\eqref{non-negative}
is then recovered by choosing $\mu_{ij} = \beta_{i}\beta_{I+j}\Pi_{ij}.$
Our proof is
variational. We construct a functional $E(\beta)$ with the property that
$\beta$ is a \emph{critical point} of $E$ --- meaning a point where $E$ has zero
derivative --- if and only if $\beta$ satisfies equation (\ref{ast}). We then
show that $E$ has exactly one critical point in the positive \emph{orthant}
$({\mathbf{R}}_{+})^{I+J}$, and give a formula for this critical point using
the Legendre transform of a related function.

The main result of this paper is the following theorem, which solves the
Choo-Siow Inverse problem.

\begin{theorem}[Existence, uniqueness, and explicit representation of a real positive solution]\label{T:main}
If all the
entries of $\Pi= (\Pi_{ij})$ are non-negative, and those of $\nu=[\,m \mid f\,]$ are
strictly positive, then precisely one solution $\beta$ of \eqref{ast} lies in
the positive orthant of ${\mathbf{R}}^{I+J}$.
%
%
Indeed, the solution $b:=(\log \beta_1,\ldots,\log \beta_{I+J})$ satisfies
$b=(DH)^{-1}(\nu)=DH^{*}(\nu)$ where $H(b)$ and $H^{*}(\nu)$ are
smooth strictly convex dual functions on ${\mathbf R}^{I+J}$ defined by
\begin{equation}\label{Teq:H}
H(b):=\frac{1}{2}\sum_{k=1}^{I+J}e^{2b_{k}} +\sum_{i =1}^{I}
\sum_{j=1}^{J} \Pi_{ij}e^{b_{i}+b_{I+j}}%
\end{equation}
and
\begin{align}
H^{\ast}(\nu) &:=   \sup_{b\in{\mathbf{R}}^{I+J}}\left\langle \nu
,b\right\rangle -H(b).\label{Teq:H star}
\end{align}
\end{theorem}
Here $\langle\cdot,\cdot\rangle$ denotes the inner product on ${\mathbf{R}}^{I+J}$.

\begin{remark}[Unpopulated types]
In case $m_i=0$ or $f_j=0$, 
we simply reformulate the problem in fewer than $I+J$ variables, corresponding only
to the populated types.  This reformulation shows the conclusions of Theorem \ref{T:main}
extend also to population vectors $\nu=[\,m \mid f\,]$ whose entries are merely non-negative,
instead of strictly positive.
\end{remark}

Since each matrix $(\mu_{ij})$ with non-negative entries solving
\eqref{equilibrium}--\eqref{non-negative} corresponds to a solution $\beta$ of (\ref{ast})
having positive amplitudes $\beta_{i} = \sqrt{\mu_{i0}}$ and $\beta_{I+j} =
\sqrt{\mu_{0j}}$, this theorem gives the sought characterization of $(\mu
_{ij})$ by $\Pi$. 
Moreover, this characterization facilitates computing variations in the marital
arrangements in response to changes in the data $(\Pi,\nu)$:

\begin{theorem}[Comparative statics] \label{T:comparative statics}
Let the unique solution to the Choo-Siow inverse problem with exogenous data $\Pi$ and $\nu$ be given by $\beta(\Pi,\nu)$.
Then (a) the percentage change of singles $\beta_k^2$ with respect to the population
parameter $\nu_\ell$ turns out to define a symmetric and positive definite matrix
\begin{equation}
r_{k\ell} := \frac{1}{\beta_{k}^2}\frac{\partial\beta_{k}^2}{\partial\nu_{\ell}};
\end{equation}
here $k,\ell \in \{1,\ldots,I+J\}$.
This \rjo{symmetric positive-}definiteness implies, among other things, the expected monotonicity $r_{kk}>0$,
the unexpected symmetry $r_{k\ell} = r_{\ell k}$,
and more subtle constraints relating
these percentage rates of change and the corresponding substitution effects such as
$|r_{k\ell}| < \sqrt{r_{kk} r_{\ell\ell}}$.

(b) Additionally we can account for the sign, and in some cases bound the magnitude,
of each entry of the matrix $R=(r_{k\ell})$. 
To avoid trivialities,  assume no column or row of $\Pi$ vanishes,  so
no observable type of individual is compelled to remain single.  Then,
\begin{equation}\label{duller statics}
r_{k\ell} <0,
 \end{equation}
 if $k \in \{1 \dots I \}$ and $\ell \in \{I+1 \dots I+J \}$ (or vice versa).
 Second, if $k,\ell \in \{1 \dots I \}$, then
\begin{equation}\label{sharper statics}
\frac{1}{2}(\beta_k^2 + \nu_k){r_{k\ell}} 
n> \delta_{k\ell}
:= \left\{\begin{array}{ll}
           0 & {\rm if}\ k \ne\ell \\
           1 & {\rm otherwise.}\ \end{array}\right.
\end{equation}
Similarly, \eqref{sharper statics} also holds if both $k,\ell \in \{I+1 \dots I+J \}$.

\end{theorem}

These qualitative comparative statics have a simple interpretation. Increased supply of
any type $k$ of man coaxes more women into marriage \eqref{duller statics}
and decreases the number of men who wish to marry. The last statement
of the theorem says that this decrease is not merely due to the fact that there are more men.
Rather, men of type $\ell \ne k$ who
would have chosen marriage under the old regime choose to be single after the shock
\eqref{sharper statics}.
The following corollary explains how the additional women are coaxed into
marriage: it shows increased competition among men leads to larger equilibrium transfers
to their female spouses, as is further explained in subsection \S\ref{S:comparative statics}.
The corollary asserts not only monotonicity of utility transferred by
men of type $i$,  but also of the percentage who choose to remain single,  as a function
of their abundance in the population.

\begin{corollary}[
Utility transferred and non-participant fraction increase with abundance]
\label{C:transfer and percent}
For all $i \le I$, $j \le J$, and $k\le I+J$, with the hypotheses and notation of
Theorem \ref{T:comparative statics},
\begin{eqnarray}
\label{transfer increases}
\frac{\partial}{\partial \nu_i} (\eta^f_{ij}-\eta^m_{ij}) &>& 0 \\
{\rm and}\quad \frac{\partial}{\partial \nu_k} \Big(\frac{\beta_k^2}{\nu_k}\Big) &>& 0.
\label{participation decreases}
\end{eqnarray}
\end{corollary}


\section{\rjo{Derived statics and further conjectures}}
\label{S:derived conjectures}

It is possible to express many quantities of interest in terms of the symmetric
positive-definite matrix $R = (r_{ij})$ from Theorem \ref{T:comparative statics},
whose entries encode the relative change in the number of type $i$ individuals
who choose not to marry in response to a fluctuation in the
total number $\nu_j$ of type $j$ individuals in the population.  For example,
in the Choo-Siow model, the number of marriages $\mu_{ij}$ of type $i$ men to type $j$
women is given by the geometric mean \eqref{equilibrium}
of the number of singles of the two given
types times the corresponding entry in the gains matrix:
$
\mu_{ij} = \Pi_{ij} (\mu_{i0}\mu_{0j})^{1/2} = \Pi_{ij} \beta_i \beta_j.
$
Since $\Pi_{ij}$ is exogenous, we immediately obtain a formula
\begin{equation}\label{relative change in (i,j) marriages}
\frac{\p \log \mu_{ij}}{\p \nu_k} = \frac{1}{2} (r_{ik} + r_{k,I+j})
\end{equation}
showing the relative change in the number of type $(i,j)$ marriages caused
by fluctuations in the total population of type $k$ individuals is just the average of
the relative changes $r_{ik}:= 2\p (\log \beta_i) / \p \nu_k$ and $r_{I+j,k}$
in the numbers of unmarrieds of the corresponding types $i$ and $j$.

We may also consider fluctuations in the number of singles of type $k$ in response
to changes in the exogenous gains parameters $\Pi_{ij}$ when the population
$\nu$ of each type of man and woman is held fixed.  In section \ref{S:gains variations},
the implicit function theorem is used to derive
\begin{equation}\label{gains derivatives 1}
\frac{\partial \beta_{k}}{\partial \Pi_{ij}}=-\beta_{i}\beta_{I+j}(\frac{\partial \beta_{k}}{\partial \nu_{i}}+\frac{\partial \beta_{k}}{\partial \nu_{I+j}})
\end{equation}
for all $i \in \{1,\ldots,I\}$, $j \in \{1,\ldots,J\}$, and $k \in \{1,\ldots,I+J\}$, or equivalently
\begin{equation}\label{gains derivatives 2}
\frac{\partial \log \beta_{k}}{\partial \Pi_{ij}}=-\frac{\mu_{ij}}{2\Pi_{ij}}({r_{ki}+r_{k,I+j}}{}).
\end{equation}
The equation (\ref{gains derivatives 1}) has an intuitive interpretation. An increase in the
total systematic gains to an $(i,j)$ marriage (produced, for example, by
an isolated increase in the value of type $j$ marriages
to type $i$ men, or an isolated decrease in the value of remaining single) has the same effect as
decreasing the supply of the men or women of the respective types by a proportionate amount,
weighted by the geometric mean of the unmarried men and women of type $i$ and~$j$.

Now Theorem \ref{T:comparative statics} shows the
the summands above to have opposite signs,  so the sign of their
sum $r_{ki}+r_{k,I+j}$ will fluctuate according to market conditions.
If $k=i$ or $k=I+j$ however,  it is natural to conjecture that this sum is positive,
in which case \eqref{relative change in (i,j) marriages} shows the
the number of type $(i,j)$ marriages $\mu_{ij}$ would be an increasing function of
the population size of type $i$ (and independently, type $j$) individuals.  Similarly
\eqref{gains derivatives 2} then asserts the number of type $i$ and $j$ singles
to be a decreasing function of the total gains $\Pi_{ij}$ for type $(i,j)$ marriages.
This conjecture amounts to requiring each positive diagonal entry $r_{kk}$ in the
matrix $R$ to dominate each negative entry $-r_{k \ell}$ in its row or column ---
a plausible strengthening (whose proof, alas, eludes us)
of the claim $r_{k\ell}^2 \le r_{kk}r_{\ell\ell}$ established in
Theorem \ref{T:comparative statics}(a).


\subsection{Summary of comparative statics}
Before turning to the proof of the theorems listed above,  let us conclude by
recapping our comparative statics.

\begin{itemize}
\item{Increasing the men of a given type increases the number of singles of all male types,
      and decreases the number of singles of all female types.
%
\item{Increasing the number of men of a given type increases the transfer they
      must pay to any woman they marry  \rjo{\eqref{transfer increases}}.}
\item{The percentage rate of change of unmarried men of type $i$ due to increases
      in women of type $j$ is equal to the percentage rate of change of unmarried
      women of type $j$ due to increases in men of type $i$.}
\item{The marital participation rate
      $\sum^{J}_{j=1} \frac{\mu_{ij}}{m_{i}}$
      decreases with an increase in own type $m_i$.}} 
%
%
\end{itemize}

\section{A New Variational Principle (Proof of Theorem \ref{T:main})}\label{S:Thm 1}

\subsection{Variational method: existence of a solution}

Consider the function $E : \mathbf{R}^{I+J} \rightarrow{\mathbf{R}}
\cup\{+\infty\}$, defined as follows:
\begin{equation}
\label{E}E(\beta):=\frac{1}{2}\sum_{k=1}^{I+J}\beta_{k}^{2} + \sum_{i = 1}^{I}
\sum_{j =1}^{J} \Pi_{ij}\beta_{i}\beta_{I+j}- \sum_{k=1}^{I}\nu_{k} \log
|\beta_{k}|.
\end{equation}
It diverges to $+\infty$ on the coordinate hyperplanes where the $\beta_{k}$
vanish, but elsewhere is smooth.

We differentiate and observe that $\beta$ is a critical point of $E$ if and
only if \eqref{ast} holds. Notice strict positivity of the components of $\nu=
[\,m \mid f\,]$ implies the corresponding component of a solution $\beta$ to
\eqref{ast} is non-vanishing, hence no solutions occur on the coordinate
hyperplanes which separate the different orthants.
In words, the critical points of $E$ are precisely those that satisfy the
system of equations we wish to show has a unique real positive root. It
therefore suffices to show that $E(\beta)$ has a unique real positive critical
point; for then (\ref{ast}) admits exactly one real positive solution. Let us
show at least one such solution exists, by showing $E(\beta) $ has at least
one critical point: namely, its minimum in the positive orthant.

\begin{claim}
[Existence of a minimum]If all the entries of $\Pi= (\Pi_{ij})$ are
non-negative, and those of $\nu=[\,m \mid f\,]$ are strictly positive, the function
$E(\beta)$ on the positive orthant defined by \eqref{E} attains its minimum value.
\end{claim}

\begin{proof}
Since $E(\beta)$ is continuous, the claim will be established if we show the
sublevel set $B_{\lambda}:= \{\beta\in({\mathbf{R}}_{+})^{I+J} \mid E(\beta)
\le\lambda\}$ is compact for each $\lambda\in{\mathbf{R}}$. Non-negativity of
$\Pi_{ij}$ combines with positivity of $\nu_{k}$, $\beta_{k}$, and the
inequality $\log\beta_{k} \le\beta_{k} -1$ to yield
\begin{align}
E(\beta)  & \ge \sum_{k=1}^{I+J} \frac{1}{2}\beta_{k}^{2} - \nu_{k} (\beta
_{k}-1)\\
& = \frac{1}{2}\sum_{k=1}^{I+J} (\beta_{k} -\nu_{k})^{2} - (\nu_{k}-1)^{2} +
1.
\end{align}
It follows that $B_{\lambda}$ is bounded away from infinity. Since $E(\beta)$
diverges to $+\infty$ on the coordinate hyperplanes, it follows that
$B_{\lambda}$ is also bounded away from the coordinate hyperplanes --- hence compactly
contained in the positive orthant.
\end{proof}

\subsection{Uniqueness, convexity, and Legendre transforms}

With this critical point characterization of the solution in mind, let us
observe for $\beta\in{\mathbf{R}}^{I+J}$ in the positive orthant, defining
$b_{k} := \log\beta_{k}$ implies $E(\beta) = H(b) - \langle\nu,b \rangle$,
where $H(b)$ is defined in \eqref{Teq:H}.
%
Since the change of variables $\beta_{k} \in{\mathbf{R}}_{+} \longmapsto b_{k}
= \log\beta_{k} \in{\mathbf{R}}$ is a diffeomorphism, it follows that critical
points of $H(b) - \langle\nu,b \rangle$ in the whole space ${\mathbf{R}}%
^{I+J}$ are in one-to-one correspondence with critical points of $E(\beta)$ in
the positive orthant.


On the other hand, $H(b)$ is manifestly convex, being a non-negative sum of
convex exponential functions of the real variables $b_{k}$; in fact $\Pi
_{ij}\geq0$ shows the Hessian $D^{2}H(b)$ dominates what it would be in case
$\Pi=0$, namely the diagonal matrix 
with positive entries
$diag[2e^{2b_{1}} ,\dots,2e^{2b_{I+J}}]$
along its diagonal. Thus $H(b)$ is strictly convex throughout
${\mathbf{R}}^{I+J}$, and $E(\beta)=H(b)-\langle\nu,b\rangle$ can admit only
one critical point $\beta$ in the positive orthant --- the minimizer whose
existence we have already shown. The solution $\beta$ to \eqref{ast} which we
seek therefore coincides with the unique point at which the maximum is attained.

This last fact means that $b$ maximizes the right-hand side of the following equation:
\begin{align}
H^{\ast}(\nu) &:=   \sup_{b\in{\mathbf{R}}^{I+J}}\left\langle \nu
,b\right\rangle -H(b)\label{H star}\\
& =\sup_{\beta\in({\mathbf{R}}_{+})^{I+J}}-E(\beta).\nonumber
\end{align}
The function $H^{*}$ defined pointwise by the above equation is the Legendre
transform or convex dual function of $H$; see
Appendix \ref{A:Legendre transform} for details.
It follows that the solution $b$ satisfies $\nu=DH(b)$. Thus $b=DH^{*}(\nu)$
by the duality of $H$ and $H^{*}$.
This provides an explicit formula for $b$ in terms of the derivative of $H^{*}$.

\section{Comparative Statics (Proof of Theorem \ref{T:comparative statics})}\label{S:Thm 2}

\subsection{Positive definiteness (a)}
Our representation of the solution in terms of the Legendre transform of the convex
function $H$ can be used to obtain information
about the derivatives of the solutions with respect to the population parameters $\nu$.

Suppose we wish to know how the number of marriages $\mu_{ij}=\Pi_{ij}\beta_{i}%
\beta_{I+j}$ of each type $(i,j)$ varies in response to slight changes in the
population vector $\nu$, assuming the gains matrix $\Pi$ remains fixed. This
is easily computed from the percentage rate of change $r_{k\ell}$ in the number $\beta
_{k}^{2}$ of unmarrieds of each type, which is given in terms of the Hessian
of either \eqref{Teq:H} or \eqref{Teq:H star} by
\begin{equation}
r_{k\ell}:=\frac{1}{\beta_{k}^{2}}\frac{\partial\beta_{k}^{2}}{\partial\nu_{\ell}}
=2D_{k\ell}^{2}H^{\ast}(\nu)=2(D^{2}H|^{-1}_{(\log \beta_1, \ldots, \log \beta_{I+J})})_{k\ell},
\quad1\leq k,\ell\leq I+J.\label{comparative statics}%
\end{equation}

To see that these equalities hold, observe that the solution $\beta$ is the point
where the maximum (\ref{Teq:H star})
is attained. The Legendre transform $H^{\ast}(\nu)$ of $H$ defined by this
maximum is manifestly convex, and its smoothness is well-known to follow from
the positive-definiteness of $D^{2}H(b)>0$ as in Lemma \ref{L:Legendre}.
Moreover $b=DH^{\ast}(DH(b))$,
whence the maximum \eqref{Teq:H star} is attained at $b=DH^{\ast}(\nu)$ and
$D^{2}H(b)^{-1}=D^{2}H^{\ast}(DH(b))=D^{2}H^{\ast}(\nu)>0$.  This positive definiteness implies the first half of the Theorem 2.

\subsection{Qualitative characterization of comparative statics (b)}
To complete our qualitative description of the substitution effects in this section,
we apply the following theorem from functional analysis to matrices
$T:\mathbf{R}^n \longrightarrow \mathbf{R}^n$.  \rjo{We define the operator norm
of such a matrix by $\displaystyle \|T\|_{op} := \max_{0 \ne v \in \mathbf{R}^n} |T(v)|/|v|$,
where $|v|=\langle v,v \rangle^{1/2}$ denotes the Euclidean norm.}

\begin{theorem}[Neumann series for the resolvent of a linear contraction]\label{neumann}
If $\|T\|_{op}<$1 \rjo{for $T:\mathbf{R}^n \longrightarrow \mathbf{R}^n$,}
the operator $(1-T)^{-1}$ exists and is equal to $\sum_{k=0}^{\infty}T^{\rjo{k}}$.
\end{theorem}

Next, we consider the matrix  $D^{2}H(b)|_{(\log \beta_1, \ldots, \log \beta_{I+J})}$,
and derive properties of its inverse, whose entries give the various values of $r_{k\ell}/2$.
Differentiating the known function $H(b)$ twice yields a positive-definite
$(I+J)\times (I+J)$ matrix which can be factored into the form
\begin{equation}\label{diagonal conjugation}
2 R^{-1} = D^2H|_{b=(\log \beta_1, \ldots, \log \beta_{I+J})} =
\Delta \begin{pmatrix}
\Delta _{I} & \Pi \\
\Pi^{T} & \Delta _{J}%
\end{pmatrix}%
\Delta
\end{equation}
where $\Delta = \diag[ e^{b_1},\ldots,e^{b_{I+J}}] = \diag[\beta]$, while $\Delta_I$
and $\Delta_J$ are $I \times I$ and $J \times J$ diagonal submatrices
whose diagonal entries are all larger than two:
\begin{eqnarray*}
(\Delta _{I})_{ii}&=&
2+\frac{1}{\beta_{i}^2}\sum_{j=1}^{J}\Pi _{ij}\beta_{i}\beta_{I+j}%
\ = 1 + \frac{\nu_i}{\beta_{i}^2},
\\
(\Delta _{J})_{jj} &=&
2+\frac{1}{\beta_{I+j}^2}\sum_{i=1}^{I}\Pi _{ij}\beta_{i}\beta_{I+j}%
\ = 1 + \frac{\nu_{I+j}}{\beta_{I+j}^2}.
\end{eqnarray*}
Here we have used the fact that the values $\beta$ are critical points and therefore satisfy the
first order conditions  \eqref{ast} 
to simplify these diagonal terms.

There are determinant  and inverse formulae for block matrices which assert \cite{Horn} that
\begin{equation}\label{det}
\det
\begin{pmatrix}
{\Delta _{I}} & \Pi \\
\Pi ^{T} & {\Delta _{J}}%
\end{pmatrix}%
=\det ({\Delta _{I}})\det ({\Delta _{J}})
\det(1-{\Delta _{I}}^{-1}\Pi {\Delta _{J}}^{-1} \Pi^T),
\end{equation}
 and
 \begin{equation}\label{inverse formula}
\begin{pmatrix}
{\Delta _{I}} & \Pi \\
\Pi ^{T} & {\Delta _{J}}%
\end{pmatrix}^{-1}
=
 \begin{pmatrix}
 ({\Delta_{I}}-\Pi{\Delta_{J}}^{-1}\Pi^{T})^{-1} & -({\Delta_{I}}-\Pi{\Delta_{J}}^{-1}\Pi^{T})^{-1}\Pi{\Delta_{J}}^{-1} \\
 -\Pi^{T}{\Delta_{I}}^{-1}({\Delta_{J}}-\Pi^{T}{\Delta^{-1}_{I}}\Pi)^{-1} & ({\Delta_{J}}-\Pi^{T}{\Delta^{-1}_{I}}\Pi)^{-1} %
 \end{pmatrix}.%
 \end{equation}
The determinant $(\ref{det})$ is positive by \eqref{diagonal conjugation} and Theorem $\ref{T:main}$.
We will now show that the eigenvalues of the
matrix $A(s)= {\Delta _{I}}^{-1}s\Pi {\Delta _{J}}^{-1}s\Pi^T$,
appearing in \eqref{det}--\eqref{inverse formula}
are bounded above by $1$ and below by $-1$ for all values of
$s\in[0,1]$. This will have implications respecting the signs of the entries of
\eqref{inverse formula},
whose $(k,\ell)^{th}$ entry is in fact equal to $\beta_k \beta_\ell r_{k\ell}/2$
hence shares the sign of the 
change \eqref{duller statics}--\eqref{sharper statics}
which we desire to estimate.
Namely, it will allow us to apply
Theorem~\ref{neumann} to block entries such as
$({\Delta_{I}}-\Pi{\Delta_{J}}^{-1}\Pi^{T})^{-1}= (1-A(1))^{-1}{\Delta_{I}}^{-1}$
in (\ref{inverse formula}).

Let $\lambda^{max}(s)$ be the largest eigenvalue of $A(s)$. Then, the smallest eigenvalue of
$(1-A(s))$ is equal to $(1-\lambda^{max}(s))$. We proceed by continuously deforming
from $s=0$ to $s=1$:
The eigenvalues of $(1-A(0))$ are equal to 1, as $A(0)$ is in fact equal to the zero matrix.
Since $\det(1-A(s))>0$ for all $s \in [0,1]$, continuity of $\lambda^{max}(s)$ and the
intermediate value theorem imply that $1-\lambda^{max}(s)>0$ for all $s$, so that
$\lambda^{max}(1)<1$. Since no row of $\Pi$ vanishes,
 $A(s)$ has positive entries whenever $s>0$. The
Perron-Frobenius theorem therefore implies that any negative eigenvalue $\lambda$ of $A(1)$
is bounded by $|\lambda|<\lambda^{max}(1)$.

Since $A$ has positive entries and $\|A\|_{op}<1$,
Theorem \ref{neumann} indicates that the entries of $(1-A)^{-1}$ are all positive ---
exceeding one on the diagonal.  But $\beta_k\beta_\ell r_{k\ell}/2$
coincides with the $(k,\ell)^{th}$ entry of
$(1-A)^{-1}\diag[\beta_1^2/(\beta_1^2+\nu_1),\ldots,\beta_I^2/(\beta_I^2+\nu_I)]$,
giving the desired inequalities \eqref{sharper statics}
whenever $k, \ell \in \{1, \dots, I\}$.
The signs of the remaining derivatives \eqref{duller statics}--\eqref{sharper statics}
may be verified by applying the same technique
to the three other submatrices present in (\ref{inverse formula}),
thus completing the proof of Theorem \ref{T:comparative statics}(a)--(b).

\section{
\rjo{Transfer utilities (Proof of Corollary \ref{C:transfer and percent} and subsequent remarks)}}
\label{S:comparative statics}


\subsection{Varying the population vectors $\nu$}
Given a specification of $\Pi$ and $\nu = [\,m_{i} \mid f_{j}\,]$, the Choo-Siow model predicts a
unique vector $\beta = [\,\mu_{i0} \mid \mu_{0j}\,]$ of unmarrieds.
Given a fixed $\Pi$ and a fixed $\beta$,
the full marriage distribution can then be uniquely recovered.
It is therefore possible to view $\mu$ as a single valued (smooth) function of $\Pi$ and $\nu$.
By Theorem \ref{T:comparative statics}, the signs of $r_{k\ell}$
are independent of $\Pi$ and $\nu$ and depend only on whether
$k \in \{1, \dots ,  I \}$, or $k \in \{I+1, \dots, I+J \}$, and likewise for $\ell$.
It is perhaps useful to visualize these comparative statics as the entries of the matrix
$D\beta$ with $D_{\ell}\beta_{k}:=\frac{\partial \beta_{k}}{\partial \nu_{\ell}}$.
Then, $D\beta$ is a block matrix that is positive in its upper-left and lower-right blocks,
and negative in its upper-right and lower-left blocks.
Schematically, 
\rjo{\eqref{duller statics}--\eqref{sharper statics} yield}
\begin{equation}\label{schematic}
D\beta = \begin{pmatrix}
+\ - \\
- \ +
\end{pmatrix}.
\end{equation}

Reverting back to the Choo-Siow notation for unmarrieds and population vectors, we have
for all $k$ and $\ell$:
\begin{equation}\label{cs}
\frac{\partial \mu_{k0}}{\partial m_{\ell}}>0, \quad
\frac{\partial \mu_{k0}}{\partial f_{\ell}}<0, \quad
\frac{\partial \mu_{0k}}{\partial f_{\ell}}>0, \quad
\frac{\partial \mu_{0k}}{\partial m_{\ell}}<0.
\end{equation}

These basic comparative statics yield qualitative information about other more complex
quantities of interest. As indicated following equation (\ref{sum}), the quantity
$\eta_{ij}^{m}+\eta_{ij}^{f}-\eta^{m}_{i0}-\eta^{f}_{0j}$ is exogenous, whereas the
first two individual summands are separately endogenous and determined within the model.
In the original formulation of this model, present in \cite{CS}, our endogenous payoff
$\eta^{m}_{ij} = \tilde \eta^{m}_{ij} - \tau_{ij}$ is separated into a systematic return
$\tilde\eta^{m}_{ij}$ presumed to be exogenous,
and a utility transfer $\tau_{ij}$ from husband to wife,
which is endogenous and set in equilibrium.
Similarly, $\eta^{f}_{ij} = \tilde \eta^{f}_{ij} + \tau_{ij}$.

In equilibrium \eqref{pre-equilibrium}, both of the following equations hold:
\begin{eqnarray}
\log(\mu_{ij})- \log(\mu_{0j})&=&\eta^{f}_{ij}-\eta^{f}_{0j}=\tilde\eta_{ij}^f+\tau_{ij}-\tilde\eta^{f}_{0j},\\
\log(\mu_{ij}) -\log(\mu_{i0})&=&\eta^{m}_{ij}-\eta^{m}_{i0}=\tilde\eta_{ij}^m-\tau_{ij}-\tilde\eta^{m}_{i0};
\end{eqnarray}
there is no utility transferred by remaining single.  Subtracting one from the other, we see that:
\begin{equation}\label{transfer}
\log(\frac{\mu_{i0}}{\mu_{0j}})=2\tau_{ij} + c_{ij},
\end{equation}
where $c_{ij}=(\tilde\eta_{ij}^f-\tilde\eta^{f}_{0j}-\tilde\eta_{ij}^m+\tilde\eta^{m}_{i0})$
is exogenous.

We denote the differentiation operator $\frac{\partial}{\partial \nu_{k}}f$ by $\dot{f}$
(suppressing the dependence on $k$). Differentiating
$c_{ij}= (\eta_{ij}^f - 2\tau_{ij}-\tilde\eta^{f}_{0j}-\eta_{ij}^m + \tilde\eta^{m}_{i0})$ and
\eqref{transfer} yields:
\begin{eqnarray*}
\frac{\partial}{\partial \nu_k}(\eta_{ij}^f -\eta_{ij}^m)
= 2 \dot{\tau_{ij}}
= \frac{\dot{\mu_{i0}}}{{\mu_{i0}}}-\frac{\dot{\mu_{0j}}}{\mu_{0j}}.
\end{eqnarray*}
The inequalities \eqref{cs} now determine the sign of $\dot{\tau_{ij}}$,
which depends on the differentiation variable $\nu_{k}$.
Since $\dot{\mu_{i0}}$ and $\dot{\mu_{0j}}$ have opposite signs, according to Theorem 2,
we find
\begin{equation}
\frac{\partial \tau_{ij}}{\partial m_{i}}>0,
\end{equation}
which means the transfer of type $i$ men to each type of spouse must increase in response
to an isolated increase in the population of men of type $i$. This is expected because an
increase in the number of type $i$ men introduces additional competition for each type of
women, due to the smearing present in the model. To decrease the number of type $i$ men
demanding marriage to a particular type of woman to a level that permits one-to-one matching
requires an increase in the transfer to crowd out some men.

While in principle the men might re-distribute so that the proportion of married men remains
the same, our next computation shows this is not the case.
We consider the marital participation rate of type $k$ individuals, or rather the
non-participation rate $s_{k}(\nu):={\beta_k^2}/{\nu_k}$,
defined as the proportion of individuals who choose not to marry.
Differentiation yields 
\begin{eqnarray*}
\frac{\partial s_k}{\partial \nu_k}
&=&   \frac{\beta_k^2}{\nu_k^2}\Big(\nu_k r_{kk} -1\Big)
\\ &>& \frac{\beta_k^2}{\nu_k^2}\Big(\frac{\nu_k - \beta_k^2}{\nu_k + \beta_k^2}\Big), 
\end{eqnarray*}
according to \eqref{sharper statics}.
But this is manifestly positive since the number $\beta_k^2$ of singles of type $k$
cannot exceed the total number $\nu_k$ of type $k$ individuals.
This means, for example, that an increase in the total population of type $k$ men increases the
percentage of type $k$ men who choose to remain unmarried,
given a fixed population of women and men of other types
(and assuming, as always, that the exogenous gains matrix $\Pi$ remains fixed).
It concludes the proof of Corollary \ref{C:transfer and percent}.

\subsection{Varying the gains data $\Pi$
\rjo{(Proof of \eqref{gains derivatives 1}-\eqref{gains derivatives 2})}}
\label{S:gains variations}

The population vector $\nu$ is one variable of interest. However the function $\beta$ also
depends on the gains parameters $\Pi_{ij}$. The complete derivative
$D_{(\nu,\Pi)}\beta = [D_\nu \beta \mid D_\Pi \beta]$
is an $(I+J) \times (IJ+I+J) $ matrix. As such there are linear dependencies among its
rows and columns. Since the matrix $D_{\nu}\beta$ is invertible, its columns are linearly
independent and form a basis of the column space. Hence, the remaining columns of
the complete derivative can be expressed using linear combinations of them. The implicit
function theorem applied to this problem turns out to yield the simple linear
relationship \eqref{gains derivatives 1}:
\begin{equation}\label{gains derivatives}
\frac{\partial \beta_{k}}{\partial \Pi_{ij}}=-\beta_{i}\beta_{I+j}(\frac{\partial \beta_{k}}{\partial \nu_{i}}+\frac{\partial \beta_{k}}{\partial \nu_{I+j}})
\end{equation}
for all $i \in \{1,\ldots,I\}$, $j \in \{1,\ldots,J\}$, and $k \in \{1,\ldots,I+J\}$.

Equilibrium \eqref{ast} coincides with vanishing of the function
$F(\beta,\nu,\Pi): \mathbf{R}^{(I+J)+(I+J)+(IJ)} \rightarrow \mathbf{R}^{I+J}$
defined by
\begin{eqnarray}
F_{i}(\nu,\Pi)=\beta_{i}^{2}+\sum_{j=1}^{J}\beta_{i}\beta_{I+j}\Pi_{ij}-\nu_{i}, && 1\leq
i\leq I\cr
F_{j}(\nu,\Pi)=\beta_{I+j}^{2}+\sum_{i=1}^{I}\beta_{i}\beta_{I+j}\Pi_{ij}-\nu_{I+j}, && 1\leq
j\leq J.
\end{eqnarray}
The implicit function theorem stipulates that if the derivative
$D_{\beta}F|_{\beta_{0},\nu_{0},\Pi_{0}}$ is invertible, there is a small neighbourhood
around $(\beta_{0},\nu_{0},\Pi_{0})$ inside which for each $(\nu,\Pi)$ there is a unique $\beta$
satisfying 
equation (\ref{ast}), and further that $\beta$ depends smoothly on $(\nu,\Pi)$.
The implicit function theorem also provides a formula for the derivative of the
implicit function
$\beta(\nu,\Pi)$. It is obtained by applying the chain-rule to $F(\beta(\nu,\Pi),\nu,\Pi)$:
\begin{equation}
[D_\nu\beta \mid D_\Pi \beta]_{\nu_{0},\Pi_{0}}=-[D_{\beta}F]^{-1}[D_\nu F \mid D_\Pi F]_{\beta_{0},\nu_{0},\Pi_{0}}.
\end{equation}
Since $\frac{\partial F_{k}}{\partial \nu_{\ell}}=-\delta_{k\ell}$,
and
$\frac{\partial F_{\ell}}{\partial \Pi_{ij}}=\beta_{i}\beta_{I+j}(\delta_{i\ell}+\delta_{I+j,\ell})$,
the first part of the preceding formula yields
$[D_\beta F]^{-1} = D_\nu \beta$, and the second part then implies \eqref{gains derivatives}.
%
%
Theorem~\ref{T:comparative statics} shows $D_\nu \beta$ is invertible, so
the hypotheses of the implicit function theorem are globally satisfied
and 
our calculations are valid.
\rjo{This concludes the proof of \eqref{gains derivatives 1}-\eqref{gains derivatives 2}.}



\appendix

\section{Derivation of the preference probabilities}\label{A:preference probability}

The random variable present in the definition of male and female utility is
the Gumbel extreme value distribution, introduced to the economics literature
by McFadden \cite{McFadden}:
\begin{definition}[Gumbel distribution]\label{ev}
A random variable $\epsilon$ is Gumbel if it has cumulative distribution function $F(\epsilon)=\exp(-\exp(-\epsilon))$.
\end{definition}
Here $Pr(\epsilon<x) = F(x)$ gives the probability that the realization of this random
variable takes a value less $x\in \mathbf{R}$.
The corresponding density function is $F'(x)=f(x)=\exp(-(x+\exp(-x))$.
The mean of $\epsilon$ is the Euler-Mascheroni constant,
which is approximately equal to $\gamma =0.57\ldots$.
Its variance is equal to $\frac{\pi^{2}}{6}$.

We now use this distribution to derive the discrete probability distribution
(\ref{probability}). 

\begin{lemma}\label{L:ev}
Suppose $\sigma>0$ and
$\eta_{ij} \in \mathbf{R}$ are constants, while for each choice of $j=0,\ldots, J,$
the $\epsilon_{ijg}$ are independent identically distributed random variables with the
Gumbel distribution. Then 
\begin{equation}\label{Boltzmann}
{\rm Pr}(\eta_{ij}+\epsilon_{ijg} = \max_{0 \le k \le J}\eta_{ik}+\epsilon_{ikg})
=  \frac{\exp(\frac{\eta_{ij}}{\sigma})}{\sum_{k=0}^{J}\exp(\frac{\eta_{ik}}{\sigma})}.
\end{equation}
\end{lemma}

\begin{proof}
It costs no generality to assume $\sigma=1$.  Then
\begin{equation}
P
:= Pr(\eta_{ij}+\epsilon_{ijg} \ge \eta_{ik}+\epsilon_{ikg} \forall k)
= \int_{-\infty}^{\infty} d\epsilon \, F'(\epsilon)\,  \Pi_{k \neq j} F(\eta_{ij}+\epsilon - \eta_{ik}).
\end{equation}
This formula follows from Bayes' rule for conditional probability, and
independence of the various random variables involved. Substituting in the explicit
formula for the Gumbel distribution from Definition \ref{ev} yields
\begin{equation}
P=
 \int_{-\infty}^{\infty} d\epsilon \exp(-(\epsilon + \exp{(-\epsilon)}) \,\Pi_{k \neq j}  \exp(-\exp(\eta_{ik}- \eta_{ij}-\epsilon )).
\end{equation}

We make a change of variables by setting $t=\exp(-\epsilon)$, so $d\epsilon$=$-{dt}/{t}$.
Evaluating the integral in the new variables yields 
\begin{eqnarray*}
P
&=&   \int_{0}^{\infty} dt \, \exp(-t)\, \Pi_{k \neq j} \exp(-t \exp(\eta_{ik}- \eta_{ij})) \\
&=& \int_{0}^{\infty}   dt \, \exp(-t \sum_{k=0}^J\exp(\eta_{ik}- \eta_{ij}))\\
&=& \frac{1}{\sum_{k=0}^{J} \exp(\eta_{ik}-\eta_{ij})} \\
&=& \frac{\exp(\eta_{ij})}{\sum_{k=0}^{J} \exp(\eta_{ik})}
\end{eqnarray*}
as desired.
%
\end{proof}

\begin{corollary}[Expected marital preferences by observed types]
Suppose a man with observable type $i$ and (unobservable) specific identity $g$ derives utility
$V^{m}_{ijg} = \eta^{m}_{ij} + \sigma \epsilon_{ijg}$ from being married to a woman
of observable type $j$, independent of her specific identity.  If $\sigma>0$,
$\eta^m_{ij} \in \mathbf{R}$ and $\epsilon_{ijg}$ are as in Lemma \ref{L:ev},
then the probability he prefers a woman of type $j$ to all other alternatives
in $\{0,1,\ldots,J\}$ is given by \eqref{probability}.
\end{corollary}

\begin{remark}[The Boltzmann / Gibbs distribution]
The probabilities which appear in \eqref{probability} and \eqref{Boltzmann} take the
form of the {\em Boltzmann} or {\em Gibbs} distributions from
statistical physics, in which the deterministic component $\eta^m_{ij}$
of the utility derived plays the role of the energy associated with marital
state $j$,  while the strength $\sigma$ of the random component plays the role of the
physical temperature.  This connection is also discussed
by Galichon and Salani\'e~\cite{GalSal} \cite{GalSal2}.
\end{remark}

\section{The Legendre Transform}\label{A:Legendre transform}

Here some well-known results pertaining to convexity and the Legendre transform are recalled.
Let $F: \mathbf{R}^{n} \rightarrow \mathbf{R}$ be a twice continuously differentiable function;
$F$ is \textit{convex} if $\text{Hess(F)}:=D^2F\geq0$, and strictly convex if the line segment
connecting any two points on the graph of $F$ lies above the graph.
The \textit{Legendre transform} or convex dual function to $F(p)$
is denoted $F^{*}()$ and defined pointwise by:
\begin{equation}
F^{*}(q)=\sup_{p \in \mathbf{R}^{n}} \{ q \cdot p-F(p)\}.
\end{equation}
Since the supremum of affine functions is convex, it is clear that $F^{*}(q)$ is a
convex function. Additionally, the following duality result is true:

\begin{lemma}[Legendre duality]\label{L:Legendre}
Let $F \in C^{2}$ be strongly convex on ${\mathbf R}^n$, meaning {\rm Hess}$(F)>0$. Then $F^{*}$
is also twice continuously differentiable. Further, if $q=DF(p)$, then p=$DF^{*}(q)$.
\end{lemma}


\begin{thebibliography}{50}
\bibitem{Becker}
Gary Becker.
\newblock A theory of marriage, Part I.
\newblock {\em Journal of Political Economy} {\bf 81} (1973) 813--846.


\bibitem{SR}
Steven Berry and Peter Reiss.
\newblock Empirical models of entry and market structure, chapter for Volume III of the
\newblock \textit{Handbook of Industrial Organization}, edited by Mark Armstrong and
Robert Porter. New York: North Holland, 2007.


\bibitem{BSV}
Loren Brandt, Aloysius Siow and Carl Vogel.
\newblock Large shocks and small changes in the marriage market for famine born cohorts in China.
\newblock University of Toronto Working Paper, 2008.


\bibitem{ChiapporiMcCannNesheim10}
Pierre-Andre  Chiappori, Robert J.\ McCann and Lars  Nesheim.
\newblock Hedonic price equilibria, stable matching, and optimal transport: equivalence, topology, and uniqueness.
\newblock {\em Economic Theory} {\bf 42} (2010) 317--354.

\bibitem{CSW}
Pierre-Andre Chiappori, Bernard Selani\'e and Yoram Weiss.
\newblock Assortative matching on the marriage market: A structural investigation.
\newblock In process.


\bibitem{CSS}
Eugene Choo, Shannon Seitz and Aloysius Siow.
\newblock Marriage matching, risk sharing, and spousal labour supplies.
\newblock University of Toronto Working Paper, 2008.

\bibitem{CS}
Eugene Choo and Aloysius Siow.
\newblock Who marries whom and why.
\newblock \textit{Journal of Political Economy} {\bf 114} (2006) 175--201.


\bibitem{Dag}
John K.~Dagsvik.
\newblock Aggregation in matching markets.
\newblock {\em International Economic Review} {\bf 41} (2000) 27--57.

\bibitem{Decker10}
Colin Decker.
\newblock {\em When do Systematic Gains Uniquely Determine the Number of Marriages
between Different Types in the Choo-Siow Marriage Matching Model?  Sufficient Conditions for a
Unique Equilibrium.}
\newblock University of Toronto MSc. Thesis, 2010.

\bibitem{Ekeland10}
Ivar Ekeland.
\newblock Existence, uniqueness and efficiency of equilibrium in hedonic markets with multidimensional types.
\newblock {\em Economic Theory} {\bf 42} (2010) 275--315.

\bibitem{FigalliKimMcCann-econ}
Alessio Figalli and Young-Heon Kim and Robert J.~McCann.
\newblock {When is multidimensional screening a convex program?}
\newblock {\em Journal of Economic Theory} {\bf 146} (2011) 454-478.

\bibitem{Fox}
Jeremy T.~Fox.
\newblock Estimating matching games with transfers.
\newblock University of Chicago working paper, 2009.

\bibitem{Gal}
Alfred Galichon.
\newblock Discussion of A. Siow's `Testing Becker's theory of positive assortative matching'.
\newblock Milton Friedman Institute, University of Chicago. February 28, 2009.

\bibitem{GalSal}
Alfred Galichon and Bernard Salani\'e.
\newblock Matching with trade-offs: Revealed preferences over competing characteristics.
\newblock Working paper, 
2009.

\bibitem{GalSal2}
Alfred Galichon and Bernard Salani\'e.
\newblock Cupid's invisible hand: social surplus and identification in matching models.
\newblock Working paper, 
2011.

\bibitem {GOZ}Neil E. Gretsky, Joseph M. Ostroy and William R. Zame. The
nonatomic assignment model. \textit{Economic Theory} {\bf 2} (1992) 103--127.

\bibitem{Horn} Roger A.~Horn and Charles R.~Johnson.
{\em Matrix Analysis}.
Cambridge: Cambridge University Press, 1985.


\bibitem{Kuber}
Felix Kuber and Karl Schmedders.
\newblock Competitive equilibria in semi-algebraic economies.
\newblock Penn Institute for Economic Research working paper, 
2007.

\bibitem{McFadden}
Daniel McFadden.
\newblock Conditional logit analysis of qualitative choice behavior,
\newblock in \textit{Frontiers in Econometrics}, edited by Paul Zarembka,
New York: Academic Press, 1974.

\bibitem{Pollack}
Robert Pollack.
\newblock Two-sex population models and classical stable population theory.
\newblock In \textit{Convergent Issues in Genetics and Demography}, edited by
Julian Adams et al., 317-33. New York, Oxford: Oxford University Press, 1990.

\bibitem{Pollard}
John H. Pollard.
\newblock Modelling the interaction between the sexes.
\newblock \textit{Mathematical and Computer Modelling} {\bf 26} (1997) 11--24.

\bibitem{Siow09p}
Aloysius Siow.
\newblock Testing Becker's theory of positive assortative matching.
\newblock Working paper 356, University of Toronto Department of Economics, 
2009.

\bibitem{Siow08}
Aloysius Siow.
\newblock How does the marriage market clear? An empirical framework.
\newblock \textit{Canadian Journal of Economics}, {\bf 41} (2008) 1121--1155.

\end{thebibliography}
\end{document}